\RequirePackage{ifpdf}
\ifpdf % We are running pdfTeX in pdf mode
\documentclass[pdftex]{sigma}
\else
\documentclass{sigma}
\fi

\begin{document}
\allowdisplaybreaks

\renewcommand{\PaperNumber}{041}

\FirstPageHeading

\ShortArticleName{Abel, Poncelet, Dirichlet and Neumann Problems}

\ArticleName{Dirichlet and Neumann Problems for String Equation,\\
Poncelet Problem and Pell--Abel Equation}

\Author{Vladimir P. BURSKII~$^\dag$ and Alexei S. ZHEDANOV~$^\ddag$}

\AuthorNameForHeading{V.P. Burskii and A.S. Zhedanov}

\Address{$^\dag$~Institute of Applied Mathematics and Mechanics
NASU, Donetsk, 83114 Ukraine}

\EmailD{\href{mailto:v30@dn.farlep.net}{v30@dn.farlep.net}}

\Address{$^\ddag$~Donetsk Institute for Physics and Technology
NASU, Donetsk, 83114 Ukraine}
\EmailD{\href{mailto:zhedanov@yahoo.com}{zhedanov@yahoo.com}}

\ArticleDates{Received November 23, 2005, in f\/inal form March
20, 2006; Published online April 12, 2006}

\Abstract{We consider conditions for uniqueness of the solution of
the Dirichlet or the Neumann problem for 2-dimensional wave
equation inside of bi-quadratic algebraic curve. We show that the
solution is non-trivial if and only if corresponding Poncelet
problem for two conics associated with the curve has periodic
trajectory and if and only if corresponding Pell--Abel equation
has a solution.}


\Classification{35L20; 14H70; 13B25}

\medskip

We will examine the Dirichlet problem for the string equation
\begin{gather}
\Phi_{xy} =0 \qquad \mbox{in}\quad \Omega, \label{SE}\\
 \Phi|_{C} =0 \qquad \mbox{on}\quad C=\partial \Omega,\label{NDP}
\end{gather}
  in a bounded
semialgebraic domain, the boundary of which is given by some
bi-quadratic algeb\-raic curve given in the Euler--Baxter form
\cite{Bax,Ves}:
\begin{gather}
F(x,y) = x^2 y^2 + 1 + a\big(x^2+y^2\big) + 2b xy =0 \label{E_Bax}
\end{gather} with real $a$ and $b$ that are
subordinated to conditions of boundedness $a>0$ and nonvanishing
$(|b|-a)^2>1$.

The Euler--Baxter curve \eqref{E_Bax} is a special case of generic
bi-quadratic curve
\begin{gather}
 \sum_{i,k=0}^2 a_{ik} x^i y^k =0
\label{gen_bi}
\end{gather} with arbitrary 9 coef\/f\/icients $a_{ik}$. In the
general situation this curve has a genus $g=1$ and hence is
birationally equivalent to  an elliptic curve \cite{GH}.

As far as we know, the John algorithm for such curve was f\/irst
considered in \cite{FrRag}, where it was shown that the
corresponding discrete dynamical system is completely integrable
and some solutions in terms of elliptic functions were presented
in special cases.

{\samepage In this paper we will  concentrate mostly on analysis
of the Euler--Baxter curve \eqref{E_Bax}, because the generic
curve \eqref{gen_bi} can be transformed to the form \eqref{E_Bax}
by elementary linear rational transformations \cite{IatRo}
\[
x \to \frac{\alpha_1 x + \beta_1}{\gamma_1 x + \delta_1}, \qquad y
\to \frac{\alpha_2 y + \beta_2}{\gamma_2 y + \delta_2}
\]
with some coef\/f\/icients $\alpha_i,\dots,\delta_i$, $i=1,2$.}

Recall that the John algorithm is def\/ined as follows
\cite{John}.

Let $\Omega$ be an arbitrary bounded domain, which is convex with
respect to characteristic directions, i.e.\ it has the boundary
$C$ intersected in at most two points by each line that is
parallel to $x$- or $y$-axes. We start from an arbitrary point
$M_1$ on $C$ and consider a vertical line passing through $M_1$.
Obviously, there are exactly two points of intersection with the
curve $C$: $M_1$ and some $M_2$. We denote $I_1$ an involution
which transform $M_1$ into $M_2$. Then, starting from $M_2$, we
consider a horizontal line passing through $M_2$. Let $M_3$ be the
second point of intersection with the curve $C$. Let $I_2$ be
corresponding involution: $I_2M_2=M_3$. We then repeat this
process, applying step-by-step involutions $I_1$ and $I_2$. Denote
$T=I_2I_1$, $T^{-1}=I_1I_2$. This transformation $T:C\to C$
produces a discrete dynamical system on $C$, i.e.\ an action of
group ${\mathbb Z}$ and each point $M\in C$ generates an orbit
$\{T^nM\,|\,n\in{\mathbb Z} \}$. This orbit can be either f\/inite or
inf\/inite. The point~$M$ corresponding to a f\/inite orbit is
called a {\it periodic} point and minimal $n$, for which $T^nM=M$,
is called {\it a period} of the point $M$. F.~John have proved
several useful assertions, among which we point out the following
one.

{\it Sufficient condition of uniqueness.} The homogeneous
Dirichlet problem has no solutions apart form the trivial one in
the space $C^2(\overline\Omega)$ if the set of periodic points of
$T$ on $C$ is f\/inite or denumerable, in particular, if there are
no periodic points on $C$.

Recall now the Poncelet algorithm \cite{Berger}. Let $A$ and $B$
be two arbitrary conics (say, two ellipses, this cases is general
and can be achieved by an appropriate projective transformation).
We start from an arbitrary point $M_1$ on the conic B and pass a
tangent to the conic $A$. This tangent intersects the conic B in
another point $M_2$. Then we pass another tangent from the point
$M_2$ to the conic $A$. We obtain then the point $M_3$ on the
conic $B$. This process can be repeated generating a set of the
points $M_1,M_2, \dots, M_n, \dots$ on the conic $B$. We denote
also $L_1,L_2,\dots$ corresponding tangential points on the conic
$A$. The famous Poncelet theorem \cite{Berger} states that if this
algorithm is periodic for $M_1$ (i.e.\ $M_1 = M_N$ for some
$N=3,4,\dots$) then this property does not depend on the choice of
initial point $M_1$ on the conic $B$. Obviously, the same property
is valid for the tangent points: $L_1=L_N$. Thus the periodicity
property for the Poncelet problem depends only on a choice of two
conics $A$ and $B$ and does not depend on a choice of the initial
point.

Our crucial observation is the following. It appears that the John
algorithm for the bi-quadratic curve $C$ is equivalent to the
Poncelet algorithm for some two conics in a two-di\-men\-sio\-nal
af\/f\/ine space.

Consider a transformation $V:(x,y)\to(u,v)$ given by the formulae
$u=x+y;$ $v=xy$. This transformation converts the bisectrix $x=y$
into the parabola $D=\{(u,v)\,|\, v=u^2/4\}$, moreover the left
halfplane $x< y$ is transformed onto the exterior $D^-$ of the
parabola $D$ and the mapping is a dif\/feomorphism, the right
halfplane $x>y$ is transformed onto the same exterior $D^-$ of the
parabola $D$, the points $(x,y)$ and $(y,x)$ have the same image.
The set of vertical lines $x=x_0$ is transformed into the set $E$
of lines $v=x_0u-x_0^2$, having the parabola $D$ as their
envelope. The set of horizontal lines $y=y_0$ is transformed also
into the same set $E$ of lines (the lines $x=c$ and $y=c$ have the
same image). The curve  \eqref{E_Bax} is transformed into the
quadric~$B$:
\begin{gather}
v^2+au^2+2(b-a)v+1=0\qquad {\rm or} \qquad
(v+b-a)^2+au^2=(b-a)^2-1.\label{Ell}
\end{gather}

Now, the John algorithm in the plane $x$, $y$ becomes the Poncelet
algorithm in the plane~$u$,~$v$, because a movement along the line
$x=c$ is converted to a movement along the tangent line to the
parabola $D$ from a f\/ixed point on the quadric $B$. The movement
from the same point of the bi-quadratic curve $C$ along the line
$y=c$ is converted to a movement along the second tangent line to
the parabola $D$ starting form the same point of the quadric $B$.
(Recall that from a generic point of a quadric $B$ there are
exactly two tangent lines to another quadric $D$.) The converse
statement is also valid.

Note that in our previous paper \cite{BurZhed2004} we proposed a
slightly dif\/ferent way to establish a connection between the
Poncelet problem and the boundary value problem \eqref{SE},
\eqref{NDP}.

This result allows one to establish some important properties of
the John algorithm for our curve $F(u,v)=0$. For example, the
Poncelet theorem says then that the periodicity property of the
John algorithm depends only on equation of the bi-quadratic curve
\eqref{E_Bax} and does not depend on the choice of initial point
on this curve.

This means that only one of two possibilities occurs:
\begin{enumerate}\vspace{-2mm}\itemsep=0pt
\item[(i)]  either all points of the curve $C$ are non-periodic.
It follows from the Poncelet problem that in this case the set
$T^nM$ is dense on $C$, it is said to be the {\it transitive}
case;

\item[(ii)] or all points of the bi-quadratic curve \eqref{E_Bax}
have the same f\/inite period $N$. \vspace{-2mm}
\end{enumerate}

Moreover, we can establish an explicit criterion for periodicity
condition for the John algorithm (see below) and to demonstrate a
relation with the famous algebraic Pell--Abel equation. This
relation comes from the well-known Cayley criterion for the
Poncelet problem. Recall that the Cayley criterion can be
formulated as follows~\cite{Berger}. Let $f(\lambda) = \det(A -
\lambda B)$ be a characteristic determinant for the one-parameter
pencil of conics $A$ and $B$ presented in the projective form. In
more details, assume that the conic $A$ has an af\/f\/ine equation
$\phi_A(x,y)=0$. We then pass to the projective co-ordinates
$x=\xi_1/\xi_0$, $y=\xi_2/\xi_0$ and present the equation of the
conic $A$ in the form
\[
\sum_{i,k=0}^2 A_{ik} \xi_i \xi_{k} =0
\]
with some $3\times 3$ matrix $A$. Similarly, the projective
equation for the conic $B$ has the form
\[
\sum_{i,k=0}^2 B_{ik} \xi_i \xi_{k} =0
\]
with some $3\times 3$ matrix $B$. We then def\/ine the polynomial
$f(\lambda) = \det(A - \lambda B)$ of the third degree. Note that
$f(\lambda)$ is a characteristic polynomial for the generalized
eigenvalue problem for two matrices $A$, $B$. Calculate the Taylor
expansion
\[
\sqrt{f(\lambda)} = c_0 + c_1 \lambda + \cdots +c_n \lambda^n +
\cdots
\]
and compute the Hankel-type determinants from these Taylor
coef\/f\/icients:
\begin{gather}
 H^{(1)}_p = \left | \begin{array}{cccc} c_3 & c_{4} & \dots &
c_{p+1}\\ c_{4}& c_{5} & \dots & c_{p+2}\\ \dots & \dots & \dots & \dots\\
c_{p+1} & c_{p+2} & \dots & c_{2p-1} \end{array} \right |, \qquad
 p=2,3,4,\dots \label{H1}
 \end{gather}
and
\begin{gather}
 H^{(2)}_p = \left | \begin{array}{cccc} c_2 & c_{3} & \dots &
c_{p+1}\\ c_{3}& c_{4} & \dots & c_{p+2}\\ \dots & \dots & \dots & \dots\\
c_{p+1} & c_{p+2} & \dots & c_{2p} \end{array} \right |, \qquad
p=1,2,3,\dots. \label{H2}
\end{gather}

 Then the Cayley criterion \cite{Berger,GH} is: {\it the
trajectory of the Poncelet problem is periodic with the period $N$
if and only if $H^{(1)}_p=0$ for $N=2p$, and $H^{(2)}_p=0$ for
$N=2p+1$.} For modern proof of the Cayley criterion see,
e.g.~\cite{GH}.
 Moreover, we have  the following observation: the
Cayley condition \eqref{H1} coincides with a solvability criterion
of the Pell--Abel equation
\begin{gather}
 A^2(\lambda)+\tilde f(\lambda)B^2(\lambda)=1\label{PA}
 \end{gather}
  with
 $\deg \tilde f=4$, $f(0)=0$ (for details concerning solvability
of the Pell--Abel equation and its relations with other problems of
mathematics see, e.g.~papers by V.A.~Malyshev
\cite{Malyshev,Malyshev2}) if one takes $ f(x)=x^4\tilde
f(x^{-1})$. Recall that in the equation \eqref{PA} the function
$f$ is known and one should f\/ind unknown polynomials $A$ and
$B$. The main problem here is to f\/ind conditions on the
coef\/f\/icient $f$ under which the equation has a solution $A$,
$B$. We have the following proposition (a detailed analysis of
this proposition will be published elsewhere):

{\it The Poncelet problem is periodic with an even period iff
corresponding the Pell--Abel equation is solvable.}

Note that the case of odd period $2N+1$ for the Poncelet problem
corresponds to the case of the same period $2N+1$ for the John
algorithm, whereas the case of even period $2N$ for the Poncelet
problem corresponds to the case of period $N$ for the John
algorithm.

Thus, in order to decide whether or not the John algorithm for the
given bi-quadratic curve~\eqref{E_Bax} is periodic we f\/irst pass
from the John algorithm to the corresponding Poncelet problem and
then apply the Cayley criterion.

Moreover, we have given the explicit solution of the John
algorithm for the bi-quadratic curve~\eqref{gen_bi} (see
\cite{BurZhed2004}):
\begin{gather} x_n = \phi(q(n - n_1)), \qquad y_n =
\psi(q(n-n_2)), \label{sol_John}
\end{gather}
 where $\phi(z)$, $\psi(z)$ are two
dif\/ferent elliptic functions of the second order with the same
periods; $n_1$, $n_2$ are some parameters depending on initial
conditions. The parameter $q$ and the periods $2\omega_1$,
$2\omega_2$ of the elliptic functions $\phi(z)$, $\psi(z)$ do not
depend on initial conditions.

These considerations allow one to write down the general solution
of the John algorithm in the form
\[
x_n = \kappa_1 \: \frac{a_1 \wp(q(n-n_1) + b_1}{c_1 \wp(q(n-n_1) +
d_1}, \qquad y_n = \kappa_2 \: \frac{a_2 \wp(q(n-n_2) + b_2}{c_2
\wp(q(n-n_2) + d_2}
\]
with some parameters $\kappa_i$, $a_i$, $b_i$, $c_i$, $d_i$,
$n_i$, $i=1,2$ \cite{BurZhed2004}.

Periodicity condition for the John algorithm is
\begin{gather} q N =
2\omega_1 m_1 + 2 \omega_2 m_2, \label{CP}
\end{gather} where $N$, $m_1$, $m_2$ are
integers.

For the special case of the Euler--Baxter biquadratic curve
\eqref{E_Bax} we have a much simpler solution
\begin{gather} x_n = \sqrt{k} \:
{\rm sn}(q(n-n_0),k), \qquad y_n = \sqrt{k} \: {\rm
sn}(q(n-n_0+1/2),k), \label{E_sol}
\end{gather} where $n_0$ is an arbitrary parameter and
parameters $q$, $k$ are easily determined from the para\-me\-ters $a$,
$b$ of the Euler--Baxter curve \cite{Bax,IatRo}.

Finally, we obtained new results concerning solutions for the
boundary value problems for the string equation in case of
periodicity of the John algorithm. For the Dirichlet problem we
have the following result: {\it Let $C$ be a nondegenerate
bi-quadratic curve \eqref{E_Bax}. The homogeneous Dirichlet
problem \eqref{SE}, \eqref{NDP} has a nontrivial solution $\Phi
\in C^2(\Omega)$ iff the John algorithm is periodic, i.e.\ the
condition \eqref{CP} is fulfilled. In this case there is an
infinite set of linearly independent smooth solutions.}

Besides, note the Dirichlet problem \eqref{NDP} has some
nontrivial solution $u $ in Sobolev spaces if\/f the Neumann
problem $u ^\prime _ {\nu_ *} \vert _ {C} = 0 $ with the conormal
$\nu_*$ for the equation \eqref{SE} has a nonconstant solution $u
$ \cite{BurBook}.

Thus, the Poncelet problem and Pell--Abel equation are closely
connected to the Neumann problem for the string equation as well.
And the condition \eqref{CP} is a criterion of existence of
nontrivial solution for each of these problems.

\LastPageEnding
\end{document}